\documentclass[twocolumn]{autart}

\usepackage{amsmath}
\usepackage{times}
\usepackage{amssymb}
\usepackage{lipsum}
\usepackage{ mathrsfs }

\usepackage{graphicx}
\usepackage{tikz}
\usepackage{pdflscape}

\makeatletter
\newcommand{\removelatexerror}{\let\@latex@error\@gobble}
\makeatother

\newtheorem{Theorem}{Theorem}

\newtheorem{Def}[Theorem]{Definition}
\newtheorem{Proposition}[Theorem]{Proposition}

\newtheorem{remark}{Remark}

\DeclareMathOperator{\Image}{Im}

\begin{document}

\begin{frontmatter}

\title{Structured Optimal Feedback in Multi-Agent Systems:  \\ A Static Output Feedback Perspective\thanksref{footnoteinfo}} 

\thanks[footnoteinfo]{Corresponding author S.~Zeng. Tel. +49 711 685 66312.}

\author[Stuttgart]{Shen Zeng}\ead{shen.zeng@ist.uni-stuttgart.de},  
\author[Stuttgart]{Frank Allg\"{o}wer}\ead{frank.allgower@ist.uni-stuttgart.de} 

\address[Stuttgart]{Institute for Systems Theory and Automatic Control, University of Stuttgart, 70550 Stuttgart, Germany}

\begin{keyword}                           
Linear quadratic optimal control, multi-agent systems, static output feedback.        
\end{keyword}                             

\begin{abstract}                          
In this paper we demonstrate how certain structured feedback gains necessarily emerge as the optimal controller gains in two linear optimal control formulations for multi-agent systems. We consider the cases of linear optimal synchronization and linear optimal centroid stabilization. In the former problem, the considered cost functional integrates squared synchronization error and input, and in the latter, the considered cost functional integrates squared sum of the states and input. Our approach is to view the structures in the feedback gains in terms of a static output feedback with suitable output matrices and to relate this fact with the optimal control formulations. We show that the two considered problems are special cases of a more general case in which the optimal feedback to a linear quadratic regulator problem with cost functionals integrating squared outputs and inputs is a static output feedback. A treatment in this light leads to a very simple and general solution which significantly generalizes a recent result for the linear optimal synchronization problem. We illustrate the general problem in a geometric light.
\end{abstract}
\end{frontmatter}

\maketitle
\thispagestyle{empty}
\pagestyle{empty}

\section{Introduction}
Structure portrays one of the key features in multi-agent systems as it allows e.g. for a cheaper implementation by reducing the required amount of information that needs to be communicated. An important class of structured couplings are diffusive couplings which require only relative information between agents. More specifically, for diffusive couplings the effective coupling terms between agents can be written in terms of differences between the states of the agents, so that there is no need for the agents to communicate absolute state information. It is sufficient that every agent determines its relative state to the other agents. Another structure that is recently gaining more attention is that of low-rank feedback gains \cite{madjidian2014distributed} of which diffusive couplings are in fact a special case. A control with a low-rank feedback structure is distinguished by a diagonal part complemented by a rank-one coordination term which corresponds to computing a simple averaging operation and then broadcasting this value to all systems uniformly \cite{madjidian2014distributed}. Such a control mechanism is very attractive as it is scalable and readily allows for a simple plug and play operation.

Diffusive couplings play a key role in the study of multi-agent systems ever since their appearance in the first works on multi-agent consensus \cite{fax2004information,olfati2004consensus}. It is thus interesting to ask whether it can be actually shown that such structures emerge necessarily as optimal feedback gains to certain optimal distributed control problems. Establishing such a result would thereby further strengthen the relevance of these structures also from an optimal control perspective. Furthermore, a result in these directions is also relevant from a practical point of view as it could justify omitting additional constraints that are introduced for the purpose of enforcing structure.

A question that goes in this direction has been first considered in the recent paper \cite{Montenbruck}. Therein an optimal control problem for the synchronization of a group of linear systems is considered with a focus on structured feedback gains. The considered optimal control problem is given by a cost functional that integrates synchronization error and quadratic input signals, and one of the results in \cite{Montenbruck} showed that in special cases, the optimal control problem results in couplings which are inherently diffusive. While previous works in this direction have already considered similar optimal control formulations \cite{cao2010optimal}, \cite{fardad2014design}, these formulations impose explicitly that the couplings be diffusive. In contrast, no such constraint is imposed in \cite{Montenbruck}, but rather it is rather shown that diffusive couplings \emph{necessarily} emerge as the solution to the optimal synchronization problem. This is also the general theme discussed in this paper.

The solution in \cite{Montenbruck} to this basic problem is based on showing that the strong (and maximal) solution of the algebraic Riccati equation in the considered linear optimal synchronization setup  has a diffusive structure. This starting point leads to a discussion of minimal positive semidefinite, strong and maximal solutions of the algebraic Riccati equation and its connections to the optimal solution in different cases; see the appendix for a summary of relevant results in linear quadratic optimal control. In special cases, the gap between the three types of solutions is closed, thus yielding the claimed result. More general cases such as when the system matrices are unstable have unfortunately escaped the scope of the results.

In this paper we demonstrate how the first result in \cite{Montenbruck} for homogeneous systems can be viewed as a special case of a linear quadratic regulator (LQR) problem in which the cost functional integrates over squared \emph{output} and input signals, and for which the resulting optimal controller is a static output feedback controller due to given circumstances. This starting point leads to a much simpler derivation which is also more general and allows to easily obtain the optimal solution. Since it is furthermore well-known that static output feedback and structured feedback gains are inherently related concepts, see e.g. \cite{syrmos1997static}, this novel approach is also appealing from a conceptual point of view. 

With this viewpoint, we also consider the optimal centroid stabilization problem and show that in our setup, the structure of low-rank feedback gains which realize a broadcast feedback emerges naturally as the optimal control. A very related problem was considered in \cite{madjidian2014distributed} where the focus is on coordinating the center of mass (centroid) of a group of homogeneous systems. Therein, the structure is essentially enforced by a constraint on the inputs of the individual agents. The focus of our results on the other hand is to show that in a similar but slightly different setup, rank-one structures necessarily emerge as the optimal solution to optimal centroid stabilization without having imposed any constraints. 

We furthermore stress that it turns out that the two considered optimal control formulations for the multi-agent systems generate what one may refer to as ``centralized'' solutions. In the optimal centroid stabilization problem, for instance, the optimal broadcast control is computed from the mean of all states of the systems. Although it might be argued from the more classical point of view of distributed control that such a control structure is undesirable, recently there has been the realization that such a control structure is in fact not impractical when the central leader only needs to perform simple computations such averaging the states of the multi-agent system. Moreover, it is one input signal computed from the aggregated (or compressed) information in terms of the mean of the systems' states that is eventually sent to all the systems. In particular, the control inputs of every system do not need to be computed individually by the central controller; they are the same for every system. Such control mechanisms are in fact steadily gaining more attention in recent years \cite{madjidian2014distributed,brockett2010control}.

The remainder of this paper is organized as follows. In Section~\ref{sec:problem_formulation} we formulate the linear optimal synchronization problem considered in \cite{Montenbruck} and discuss our treatment in the output-based framework. Our approach and result is illustrated in a simple example which also initiates a comparing discussion between the novel solution and the previous solution. In Section~\ref{sec:general_output} we discuss the general underlying principle that  is used in our treatment of the two specific problems, which is of separate interest. The general principle is illustrated from a geometric point of view. Using this same principle, we show in Section~\ref{sec:centroid} how broadcast feedback is the optimal feedback gain for optimal centroid stabilization. Lastly, Section~\ref{sec:conclusions} concludes the paper.

\section{Necessity of diffusive couplings in a linear optimal synchronization problem}
\label{sec:problem_formulation}

\subsection{Problem formulation}
We consider a group of $N$ homogeneous linear systems
\begin{align*}
 \dot{x}_i = A x_i + B u_i
\end{align*}
where $A \in \mathbb R^{n \times n}$, $B \in \mathbb R^{n \times p}$, and $(A,B)$ is controllable. Synchronization of the agents in the group refers to the situation when
\begin{align*}
   \lim_{t \to \infty}  \| x_i - x_j \| = 0
\end{align*}
for all $i,j  \in \{1, \dots, N \}$. This can also be described by the fact that the solution of the stacked system
\begin{align}
   \dot{x} = (I_N \otimes A) x + (I_N \otimes B) u
\label{eq:stacked}
\end{align}
approaches the synchronization subspace
\begin{align*}
  \mathscr{S} := \Image  ( \mathbf{1}_N \otimes I_n),
\end{align*}
where $\mathbf{1}_N$ denotes the all ones vector of length $N$.
The synchronization error $e$ is then defined as the orthogonal projection of $x$ onto the orthogonal complement of $\mathscr{S}$ which is called the asynchronous subspace $\mathscr{A}$. By introducing $P_{\mathscr{A}}$ as the orthogonal projection on $\mathscr{A}$, the synchronization error is given by $e = P_{\mathscr{A}} x$.

The cost functional considered in \cite{Montenbruck} is accordingly given by
\begin{align}
   J = \int_{0}^\infty  (x^{\top} P_{\mathscr{A}}^{\top} Q P_{\mathscr{A}} x + u^{\top} R u ) \, \text{d}t,
\label{eq:cost_montenbruck}
\end{align}
where $Q, R > 0$. The linear optimal synchronization problem then amounts to finding the control $u$ that minimizes the cost functional subject to \eqref{eq:stacked}. Similarly as in \cite{Montenbruck} we allow for $u$ for which the resulting state trajectory of \eqref{eq:stacked} does not converge to the origin. This is to allow for nontrivial solutions, such as periodic orbits on $\mathscr S$. 

The first result in \cite{Montenbruck} for the homogeneous case states that for specific classes of system matrices $A$, any optimal control is necessarily diffusive. Thereby the abstract definition of a diffusive coupling $u = Kx$ is via the following definition.
\begin{Def}[cf. \cite{Montenbruck}]
A matrix $K \in \mathbb R^{Np \times Nn}$ that satisfies
\begin{align}
  K  (\mathbf{1}_N \otimes I_n) = 0
\label{eq:structural_constraint}
\end{align}
 is said to be \emph{diffusive}.
\end{Def}

Next we illustrate a new way of viewing this as an optimal control problem involving static output feedback, and based on this approach, give a simpler proof which furthermore results in more general and more detailed insights in the case of homogeneous systems.

\subsection{Diffusive couplings as a static output feedback}
\label{sec:main}
The key idea of our approach is to consider the problem in a static output feedback framework. A first clue that static output feedback could be relevant is the claimed (necessary) structural constraint \eqref{eq:structural_constraint} on the feedback gain, cf. \cite{Montenbruck}. It is well-known that static output feedback can be viewed as a state feedback that is subject to \emph{structural constraints} \cite{syrmos1997static}. These are described by the condition $KY = 0$, where $Y$ is an orthonormal basis of the kernel of the output matrix $C$. By \eqref{eq:structural_constraint}, this suggests to consider a matrix $C$ such that $\ker C = \mathscr S$. Thus, we consider $C$ whose rows are given by the orthonormal basis of $\mathscr A$. Then the vector $Cx$ consists of the components of a vector $x$ in the directions of the orthonormal basis of $\mathscr A$ which is, similarly as $P_{\mathscr A} x$, a reasonable measure for synchronization error. In the following we introduce our notation for the orthonormal basis of $\mathscr A$ and then proceed towards formulating the considered optimal synchronization problem in an output-based approach.

Let $\mathcal L$ denote the graph Laplacian of the complete graph $K_N$. Then
$
  \ker \mathcal L = \Image \mathbf{1}_N
$
is the synchronous subspace for $N$ scalar systems, which is denoted by $\mathscr S'$. The orthogonal complement of $\mathscr S'$ is called the asynchronous subspace for $N$ scalar systems, which is denoted by $\mathscr A'$. Now let $\Gamma_1, \dots, \Gamma_{N-1} \in \mathbb R^N$ denote an orthonormal basis of $\mathscr A'$ and introduce the matrix
\begin{align*}
   \Gamma = \begin{pmatrix} \Gamma_1 & \dots & \Gamma_{N-1} \end{pmatrix}.
\end{align*}
As motivated in the beginning of this section, we consider the fictitious output
\begin{align*}
    y = ( \Gamma^{\top} \otimes I_n) x =: Cx.
\end{align*}
Note that the rows of $C$ form an orthonormal basis of $\mathscr{A}$. Given this output, we then consider for some $\tilde{Q}>0$ the cost functional
\begin{align*}
 J = \int_{0}^{\infty} (y^{\top} \tilde{Q}  y + u^{\top} R u ) \, \text{d}t.
\end{align*}

We will show that the input signal that minimizes the above cost functional is a static output feedback $u= Ky$, which yields the claimed diffusive structure by construction. To this end, a crucial step is to consider
\begin{align*}
  \dot{y} = (\Gamma^{\top} \otimes I_n) \dot{x} = (\Gamma^{\top} \otimes I_n) \left(  (I_N \otimes A) x + (I_N \otimes B) u  \right).
\end{align*}
Due to the block diagonal structure of the matrices $(I_N \otimes A)$ and $(I_N \otimes B)$, and the compatible block structure of $C$,
\begin{align*}
  \dot{y} = (I_{N-1} \otimes A) (\Gamma^{\top} \otimes I_n) x +  (I_{N-1} \otimes B) (\Gamma^{\top} \otimes I_p)  u,
\end{align*}
i.e. the matrices $(\Gamma^{\top} \otimes I_n)$ and $(\Gamma^{\top} \otimes I_p)$ were pushed through, while at the same time the dimensions of the identity matrices were reduced by one.

This can also be seen from considering the equation
\begin{align*}
   (\Gamma^{\top} \otimes I_n)   (I_N \otimes A) &= (\Gamma^{\top} \otimes A) \\
							   &= (I_{N-1} \otimes A) (\Gamma^{\top} \otimes I_n),
\end{align*}
and similarly for the second term involving $B$.

Since we defined $C =  (\Gamma^{\top} \otimes I_n)$, with $y = Cx$ the problem that we considered initially reduces to the following optimal control problem
\begin{align}
\begin{split}
    &\underset{u(\cdot)}{\text{minimize }}  \;  \int_{0}^{\infty}   (y^{\top} \tilde{Q} y + u^{\top} R u ) \, \text{d}t    \\
    &\text{subject to }  \;\; \dot{y} = (I_{N-1} \otimes A) y + (\Gamma^{\top} \otimes B) u.
\end{split}
\label{eq:lqr_output}
\end{align}
Note that $((I_{N-1} \otimes A), (I_{N-1} \otimes B) (\Gamma^{\top} \otimes I_p))$ is controllable since $((I_{N-1} \otimes A),(I_{N-1} \otimes B))$ is, and since $ (\Gamma^{\top} \otimes I_p)$ is full row rank. Thus, if we have $\tilde{Q}, R > 0$, by standard LQR results, the optimal input $u(\cdot)$ to the problem \eqref{eq:lqr_output} must be of the form $u = Ky$. The remarkable fact that lead to this conclusion is that the dynamics of the output $y$ which appears in the cost functional can be written down in a closed way, i.e. the dynamics of $y$ does only depend on $y$ and $u$.

\begin{remark}
\label{remark:cqc}
To recover the direct connection to the functional \eqref{eq:cost_montenbruck}, note that since $P_{\mathscr A} = C^{\top} C$, choosing $\tilde{Q} = C Q C^{\top}$ does the job. Furthermore, if $Q$ is positive definite, then $\tilde{Q}$ is positive definite as well. In fact, with the restriction that $QC^{\top} \ne 0$, we could even allow for the weaker assumption that $Q$ in \eqref{eq:cost_montenbruck} is positive semidefinite.
\end{remark}

In particular, we can conclude with the following result.
\begin{Theorem}
\label{thm:diffusive_output}
Let $(A,B)$ be controllable, and let $Q,R > 0$. Then the solution to the optimal synchronization problem
\begin{align}
\begin{split}
      &\underset{u(\cdot)}{\text{minimize }}  \;  \int_{0}^\infty  (x^{\top} P_{\mathscr{A}}^{\top} Q P_{\mathscr{A}} x + u^{\top} R u ) \, \text{d}t  \\
      &\text{subject to }  \;\; \dot{x} = (I_N  \otimes A) x + (I_N \otimes B) u
\end{split}
\label{eq:lqr_sync_montenbruck}
\end{align}
is given by $u=Kx$ where $K$ is \emph{diffusive}.
\end{Theorem}

\begin{pf}
First of all, by Remark~\ref{remark:cqc}, we can consider
\begin{align}
  J =  \int_{0}^{\infty}   (y^{\top} \tilde{Q} y + u^{\top} R u ) \, \text{d}t
\label{eq:output_LQR}
\end{align}
instead of the cost functional in \eqref{eq:lqr_sync_montenbruck} by simply identifying $\tilde{Q} = CQC^{\top}$ and $y = Cx$. 

In the following we first argue that the optimal control to \eqref{eq:lqr_output} is indeed also the optimal control for \eqref{eq:lqr_sync_montenbruck}. The claim will then follow from  standard LQR results for the reduced system \eqref{eq:lqr_output} which immediately yields that the optimal control is of the form $u=Ky$.

Let $\tilde{x}_0 \in \mathbb R^{Nn}$ and an input $\tilde{u}: [0, \infty) \to \mathbb R^{Np}$ be given and let $\tilde{x}: [0,\infty) \to \mathbb R^{Nn}$ denote the solution to
\begin{align*}
\dot{x} = (I_N  \otimes A) x + (I_N \otimes B) u, \hspace{0.5cm}  x(0) =\tilde{x}_0.
\end{align*}
Given the solution $\tilde{x}$, we compute $\tilde{y}=C\tilde{x}$ and insert the fictitious output $\tilde{y}$ and the input $\tilde{u}$ into the cost functional \eqref{eq:output_LQR}. We denote the resulting value of the cost functional by $\tilde{J}$ provided that it is finite. Differentiating the fictitious output on the other hand yields a linear system of the form
\begin{align}
  \dot{y} = (I_{N-1} \otimes A) y +  (\Gamma^{\top} \otimes B) u.
\label{eq:y_linear_system}
\end{align}
By basic existence and uniqueness arguments, the solution of the above linear system with $y(0) = C \tilde{x}_0$ and $u = \tilde{u}$ is given by $\tilde{y}$. In other words, inserting the solution of the linear system \eqref{eq:y_linear_system} with initial condition $y(0) = C \tilde{x}_0$ and $u = \tilde{u}$ into the cost functional \eqref{eq:output_LQR} leads to the exact same value $\tilde{J}$. In particular, the optimal control of \eqref{eq:lqr_output} and the optimal control of \eqref{eq:lqr_sync_montenbruck}, provided that they exist, are equal.

This justifies to solely consider \eqref{eq:lqr_output} which has the nice feature that it is given in the form of a standard LQR problem where only $x$ is replaced by the variable $y$. We have that $\tilde{Q}, R > 0$ and that $((I_{N-1} \otimes A), (I_{N-1} \otimes B) (\Gamma^{\top} \otimes I_p))$ is controllable. Therefore, an optimal control for  \eqref{eq:y_linear_system} exists, and furthermore the optimal control is a static feedback in the variable $y$, i.e. $u=\tilde{K}y$. By our foregoing arguments, $u=\tilde{K}y = \tilde{K}Cx$ is also the optimal control for \eqref{eq:lqr_sync_montenbruck}. Furthermore $K = \tilde{K} C$ is a diffusive gain by the definition of $C$. \hfill \phantom{.} \qed
\end{pf}

With Theorem~\ref{thm:diffusive_output} we showed that the solution to the optimal synchronization problem is necessarily a diffusive coupling, which is an interesting conceptual result. In distributed control, it is furthermore desirable that the computation of the couplings can be efficiently solved and in particular is scalable. In the following we show that this is the case when the weights of the linear optimal synchronization problem \eqref{eq:lqr_output} are chosen as $\tilde{Q} = I_{N-1} \otimes V$ and $R = I_{N-1} \otimes W$ with $V, W > 0$. With this choice of weights the considered cost functional is the sum of the individual cost functionals of the agents, where the weights are homogeneous.

\begin{Proposition}
\label{prop:scalability}
Let $(A,B)$ be controllable, $\tilde{Q} = I_{N-1} \otimes V$ and $R = I_{N-1} \otimes W$ with $V, W > 0$. Then the solution to the optimal synchronization problem is a homogeneous all-to-all diffusive coupling
\begin{align*}
  u_i = \frac{1}{N} \sum_{j=1}^N   W^{-1} B^{\top} Y (x_j - x_i),
\end{align*}
where $Y$ is the unique positive definite solution to the ARE
\begin{align*}
   Y BW^{-1}B^{\top} Y - YA - A^{\top} Y - V = 0.
\end{align*}
\label{prop:diffusive_and_scalable}
\end{Proposition}

For the proof of Proposition~\ref{prop:diffusive_and_scalable} we refer to Appendix~\ref{appendix:proof}.  We point out that such an analysis was also considered in \cite{Montenbruck}, in a slightly different framework.

To conclude, we have seen that considering an output in the cost functional does have applications to e.g. synchronization problems where it is desired that the synchronization error approaches zero and not necessarily the state, and we have presented the linear optimal synchronization problem in this more general output-based framework. In the following we illustrate this in a simple example which also initiates a comparing discussion of this novel approach and the existing approach. Afterwards, we discuss our result for the optimal synchronization problem as a special case of a more general phenomenon in which the optimal control minimizing a cost functional that integrates outputs and inputs is necessarily given by static output feedback.

\subsection{Illustrative examples and comparing discussion}
\label{sec:example}
In this subsection, we illustrate our approach in a simple example taken from \cite{Montenbruck} where it was used to demonstrate some delicate points in the case that the system matrix has unstable eigenvalues. This is in contrast to our approach, in which the example has a simple and definite optimal solution which is furthermore quite easy to obtain.

We consider a situation with two scalar agents $\dot{x}_i = x_i + u_i$ and where the matrix $Q$ and $R$ are the $2 \times 2$-identity matrices. The synchronous subspace here is given by the span of vector $\big( 1 \; 1 \big)^{\top}$. Thus the vector $C^{\top} = \frac{1}{\sqrt{2}} \big(-1 \; 1 \big)^{\top}$ forms an orthonormal basis of the asynchronous subspace $\mathscr A'$. With $\tilde{Q} = C Q C^{\top}$, this results in the cost functional
\begin{align*}
  J = \int_{0}^\infty (y^2 + u^{\top} u ) \, \text{d}t.
\end{align*}
Furthermore, the dynamics of the output can be written out in a closed form, since
\begin{align*}
   \dot{y} = C \Bigg( \begin{pmatrix} 1 & \\ & 1 \end{pmatrix} x +   \begin{pmatrix} 1 & \\ & 1 \end{pmatrix}u \Bigg)
		 = y + \frac{1}{\sqrt{2}} \begin{pmatrix} -1 & 1 \end{pmatrix} \begin{pmatrix} u_1 \\ u_2 \end{pmatrix}.
\end{align*}

The resulting optimization problem is thus given by
\begin{align*}
      &\underset{u(\cdot)}{\text{minimize }}  \;  \int_{0}^\infty (y^2 + u^{\top} u ) \, \text{d}t  \\
      &\text{subject to }  \;\;  \dot{y} = y + \frac{1}{\sqrt{2}} \begin{pmatrix} -1 & 1 \end{pmatrix} \begin{pmatrix} u_1 \\ u_2 \end{pmatrix}
\end{align*}
which satisfies all the requirements of the standard LQR problem.

The solution of the algebraic Riccati equation is given by $P = 2.4142$, and the resulting feedback is given by
\begin{align}
\begin{split}
  u &= \begin{pmatrix} \phantom{-}1.7071 \\ -1.7071 \end{pmatrix} y = \begin{pmatrix} \phantom{-}1.7071 \\ -1.7071 \end{pmatrix} \frac{1}{\sqrt{2}} \begin{pmatrix} -1 & 1 \end{pmatrix} x \\
      &= \begin{pmatrix} -1.2071 & \phantom{-}1.2071 \\ \phantom{-}1.2071  & -1.2071 \end{pmatrix} x,
\end{split}
\label{eq:numerical_solution}
\end{align}
which is diffusive since we can write
\begin{align*}
  u_1 = 1.2071(x_2 - x_1), \hspace{0.5cm} u_2 = 1.2071  (x_1 - x_2).
\end{align*}
The resulting closed loop system for our diffusive couplings is given by
\begin{align*}
  \begin{pmatrix} 1  & 0 \\ 0 & 1 \end{pmatrix} +   \begin{pmatrix} -1.2071 & \phantom{-}1.2071 \\ \phantom{-}1.2071  & -1.2071 \end{pmatrix} = \begin{pmatrix} -0.2071 & \phantom{-}1.2071 \\ \phantom{-}1.2071 & -0.2071 \end{pmatrix}
\end{align*}
which has eigenvalues $-1.4142$ and $1$. The stable eigenvalue is associated to the eigenvector $\frac{1}{\sqrt{2}}\big(-1 \; 1\big)^{\top}$ which describes the asynchronous mode of the two scalar agents. The unstable eigenvalue on the other hand corresponds to the synchronous mode, which does not add to the cost of our considered cost functional.

\begin{figure}[h!]
\vskip 0.12cm
  \centering
\hspace{-0.2cm}
   \includegraphics[width=0.245 \textwidth]{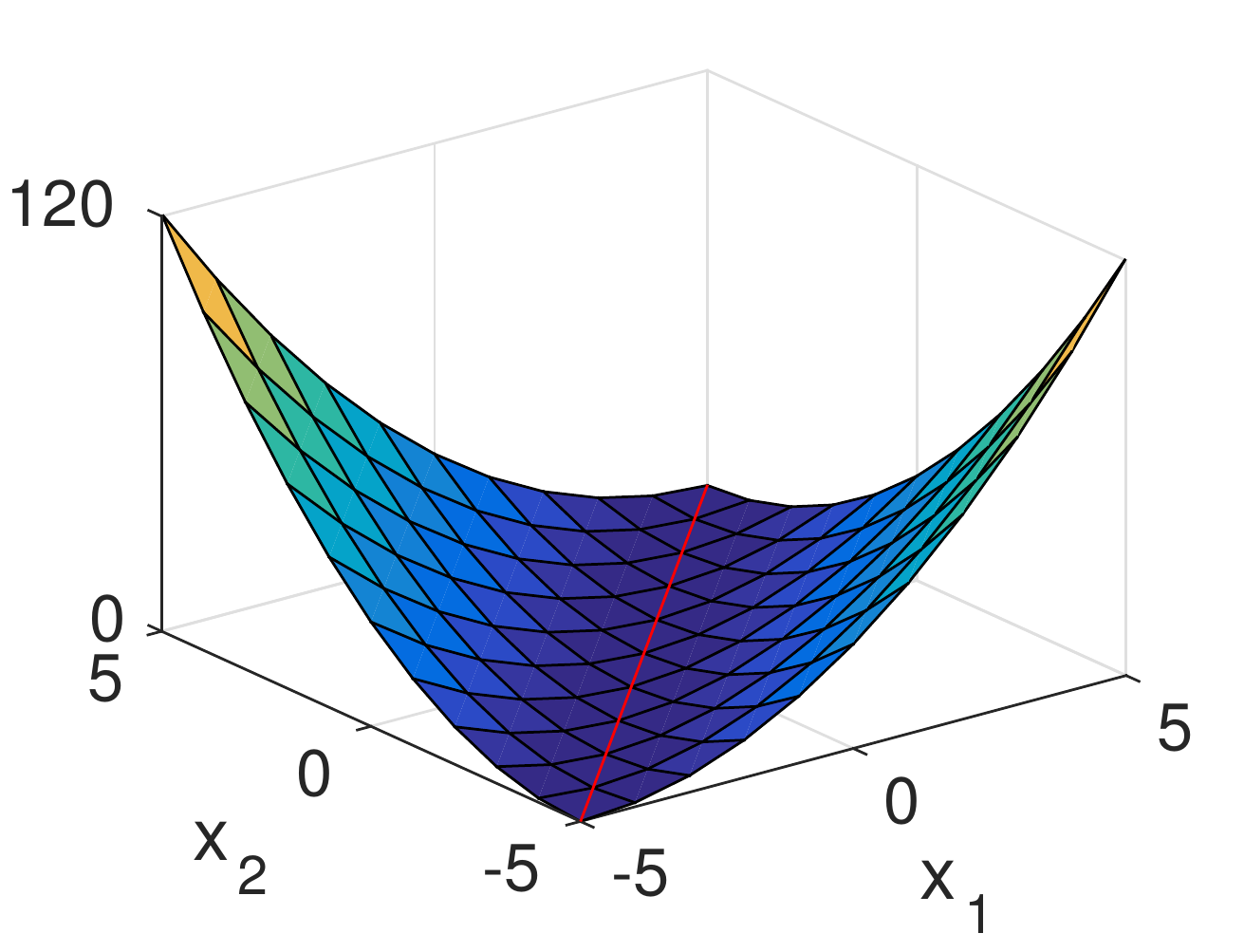}   \hspace{-0.3cm}
   \includegraphics[width=0.245 \textwidth]{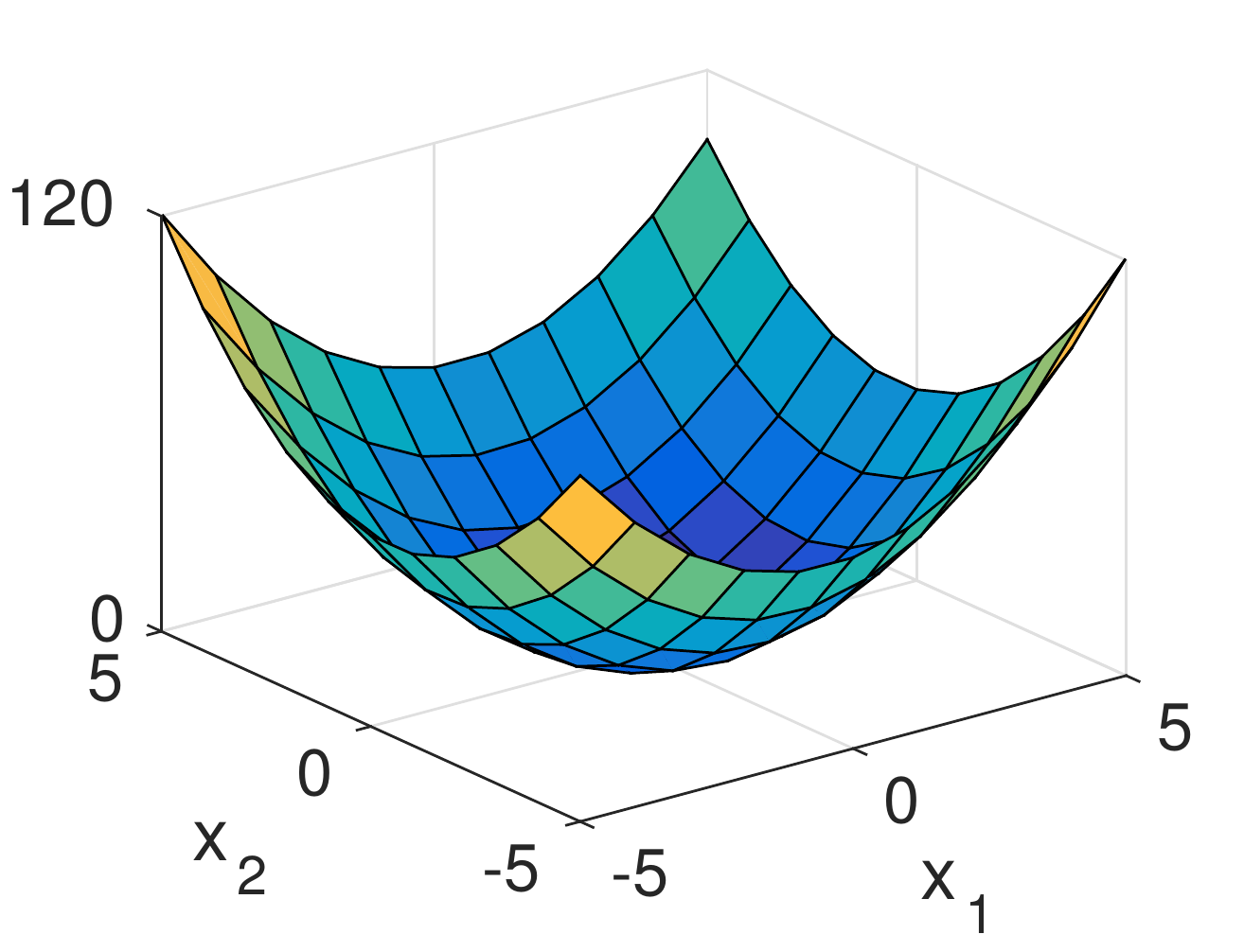} \hspace{-0.3cm}
  \caption{Top: Plot of the function $x_0 \mapsto J^{\star}(x_0) = x_0^{\top} C^{\top} P C x_0$. The cost is zero along the synchronization subspace (red line). Bottom: Plot of the function $x_0 \mapsto x_0^{\top} X_s x_0$ using the strong solution of the algebraic Riccati equation from \cite{Montenbruck}.}
 \label{fig:cost_quadratic_form}
\vskip 0.15cm
\end{figure}

Since the approach in \cite{Montenbruck} was based on structured strong solutions of the algebraic Riccati equation in which the spectral conditions of $A$ were already crucial, results for cases with $A$ having unstable eigenvalues could not be formulated. 
But it was shown for this example that the strong solution of the algebraic Riccati equation 
\begin{align*}
   X_s = \begin{pmatrix} \phantom{-}2.2071 & -0.2071 \\ -0.2071 & \phantom{-}2.2071 \end{pmatrix}
\end{align*}
is not diffusive, and that the resulting feedback $u' =  - X_s x$ in \cite{Montenbruck} is hence also not diffusive. To compare the resulting cost of the two approaches to the linear optimal synchronization problem, we compute
\begin{align}
   X^{\star} := C^{\top}    P C = \begin{pmatrix} \phantom{-}1.2071 & -1.2071 \\ -1.2071 & \phantom{-}1.2071 \end{pmatrix},
\label{eq:cost_quadratic_form}
\end{align}
which follows from the fact that for a given initial state $x_0  \in \mathbb R^2$ the resulting cost is given by $Py_0^2$ where $y_0=Cx_0$.
In Figure~\ref{fig:cost_quadratic_form} we compare the cost given by the quadratic form induced by \eqref{eq:cost_quadratic_form} to the cost given by the quadratic form that is induced by the strong solution $X_s$. In one direction the costs seem to behave similarly. In the other direction, the  feedback $u' = -X_s x$ from \cite{Montenbruck} results in non-zero costs along the synchronous subspace, which is in contrast to our solution. 

Indeed, the difference of the two matrices $X^{\star}$ and $X_s$ is
\begin{align}
    X^{\star} - X_s =   \begin{pmatrix}   -1 &    -1\\      -1  &  -1 \end{pmatrix},
\label{eq:difference_negative_definite}
\end{align}
which is negative semidefinite, i.e. $X^{\star} \le X_s$. 

In the following we try to clarify the delicate points in the approach in \cite{Montenbruck} that considers the optimal synchronization problem  from the beginning as a generic linear quadratic regulator problem
\begin{align}
\begin{split}
      &\underset{u(\cdot)}{\text{minimize }}  \;  \int_{0}^\infty  (x^{\top} \check{Q} x + u^{\top} R u ) \, \text{d}t  \\
      &\text{subject to }  \;\; \dot{x} = (I_N  \otimes A) x + (I_N \otimes B) u
\end{split}
\label{eq:LQR_general}
\end{align}
with $\check{Q} = P_{\mathscr{A}}^{\top} Q P_{\mathscr{A}}$. First of all we recall that the strong solution $X_s$ given in \cite{Montenbruck} is neither diffusive, nor optimal (see our optimal solution). In the two cases of only stable poles and only purely imaginary poles respectively, the strategy in \cite{Montenbruck} was to show that the strong solution of the algebraic Riccati equation is structured. Under the spectral conditions of $A$, the strong solution is also the maximal solution which indeed yields the optimal solution. In the case that the matrix $A$ has unstable eigenvalues, this example served as an example where the strong solution is not diffusive. 

However, the relation of this strong solution to the linear optimal synchronization problem was left open. Through our solution, we are now able to conclude that the feedback of \cite{Montenbruck} that is obtained by the strong solution (i.e. the maximal solution) is not optimal and therefore actually irrelevant.

The precise reason is that we can only infer the optimal control from the maximal solution in the special case that the system matrix has eigenvalues in the closed left half-plane which is of course not the case in this example. The solution of the algebraic Riccati equation that is more appropriate to consider in the general case is the smallest positive semidefinite solution $X_{-}$ which yields always the optimal feedback. Therefore, we conclude that the delicate situation is attributed to a gap $X_{-} \ne X_{+}$. The reason why such gap occurs, is answered by Theorem~\ref{thm:all_three_solutions} in the appendix which provides a necessary and sufficient condition.

To conclude, our novel approach completely circumvents the discussion of (structured) strong solutions of the algebraic Riccati equation for \eqref{eq:LQR_general} with $\check{Q} = P_{\mathscr A}^{\top} Q P_{\mathscr A}$ and thus also the obstructions in different cases such as when the system matrix has unstable eigenvalues. We proved that our approach does indeed yield the optimal solution. By inherently formulating the problem as a static output feedback problem, we are able to obtain \emph{the optimal} solution directly. Moreover, it is also very easy to compute the optimal solutions numerically.

\section{The general problem and its geometric interpretation}
\label{sec:general_output}
While it was not apparant in the beginning, it turned out that the key to easily solve the optimal synchronization problem includes a simple fact which can be restated in a general linear systems framework as the existence of a matrix $\tilde{A}$ so that 
$
CA = \tilde{A}C.
$
This is relevant for deriving
\begin{align*}
 \dot{y} = C\dot{x} &= C(Ax+Bu) \\
	&= \tilde{A}Cx + CBu = \tilde{A}y + CBu.
\end{align*}
Furthermore, due to the structure of the problem we could show that the resulting system for the output $(\tilde{A},CB)$ is controllable. Therefore, by the same arguments of the foregoing proof, minimizing \eqref{eq:output_LQR} only requires a static output feedback $u=Ky$ for the original system $\dot{x} = Ax + Bu$. 

The more general underlying principle can be then stated as follows.
\begin{Proposition}
Let $A,B,C$ be given such that there exists a matrix $\tilde{A}$ with $CA=\tilde{A}C$ and $(\tilde{A},CB)$ controllable. Then the solution to the linear quadratic regulator problem
\begin{align*}
      &\underset{u(\cdot)}{\text{minimize }}  \;  \int_{0}^\infty (y^{\top}Qy + u^{\top}Ru ) \, \text{d}t  \\
      &\text{subject to }  \;\;  \dot{x} = Ax + Bu,
\end{align*}
with $Q,R>0$, is a static output feedback, i.e. $u = Ky$.
\end{Proposition}

A proof of this result is completely along the lines of the proof of Theorem~\ref{thm:diffusive_output}.

For a minimal example that illustrates this general idea, as well as one of its geometric interpretations, we consider the following LQR problem
\begin{align*}
\begin{split}
      &\underset{u(\cdot)}{\text{minimize }}  \;  \int_{0}^\infty  (y^2 + u^2 ) \, \text{d}t  \\
      &\text{subject to }  \;\; \dot{x} = \begin{pmatrix} a_{11} & a_{12} \\ 0 & a_{22} \end{pmatrix} x + \begin{pmatrix} b_1 \\ b_2 \end{pmatrix} u \\
      &\hspace{1.73cm}		y = \begin{pmatrix} 0 & 1 \end{pmatrix} x.
\end{split}
\end{align*}
Again, the key feature of this problem formulation is that the cost functional depends only on the input and the second state $y=x_2$, and furthermore, that the second state is completely dynamically decoupled from the first state. Therefore, it is intuitively clear that there is only the need to consider the second state, and that in particular a feedback mechanism for the optimal control only needs to take information of the second state into account. Indeed, differentiating the output, we obtain
$$
  \dot{y} = a_{22} y + b_2 u,
$$
and thus if $b_2 \ne 0$, then an optimal control will exist and will be a (scalar) static output feedback, cf. the proof of Theorem~\ref{thm:diffusive_output} for the detailed argument.

A more general systems theoretic interpretation of the existence of a matrix $\tilde{A}$ so that $CA = \tilde{A}C$ is the invariance of $\ker C$ under the mapping $x \mapsto Ax$ (see also \cite{wonham1973tracking}), and equivalently the invariance of $\ker C$ under $\dot{x}=Ax$. This is because for any $x \in \ker C$, we can derive
\begin{align*}
  C(Ax) = CAx = \tilde{A}Cx = 0,
\end{align*}
i.e. $x \in \ker C \Rightarrow Ax \in\ker C$. To illustrate the geometric interpretation of this, we consider the drift term $CAx$ of the dynamics 
$
\dot{y} = CAx + CBu
$
which describes the effect that the vector field $\dot{x}=Ax$ has on the output. Let $x'$ be arbitrary and consider the shifted vector $x'' = x' + h$ where $h \in \ker C$. Then the drift in the output at the shifted point is
\begin{align*}
CAx'' = CA(x'+h) = CAx' + C(Ah).
\end{align*}
Since $Ah \in \ker C$, we find that $CAx' = CAx''$. Thus the change in the output $y$ by virtue of the drift and input that will be experienced at any two points $x'$ and $x''$ related by $Cx'= y = Cx''$ is the same. 

\begin{figure}[h!]
\centering
 \includegraphics[width=0.36 \textwidth]{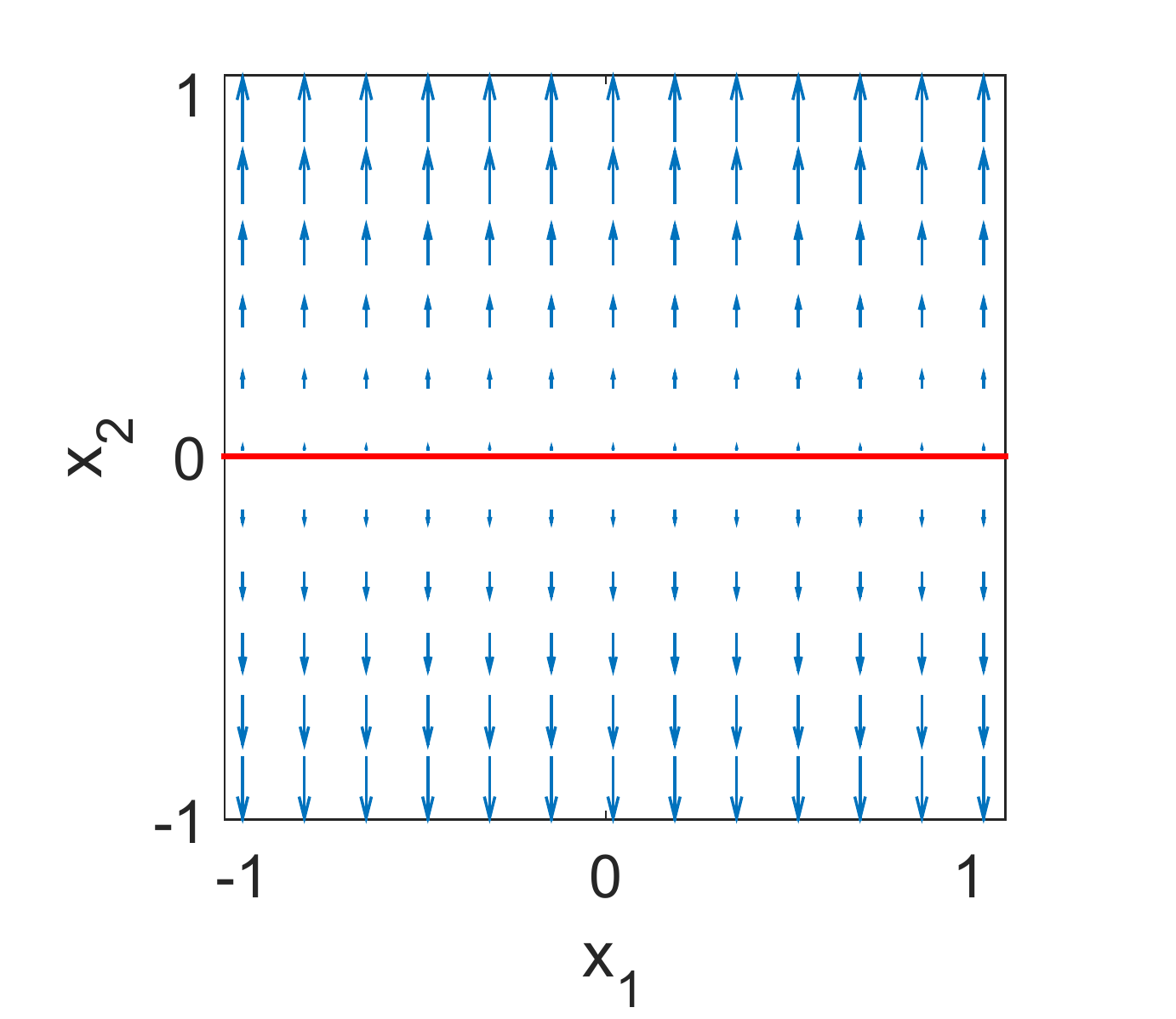}
\vspace{-0.16cm}
\caption{The strength of the vector field in the orthogonal direction of $\ker C$ is the same along affine subspaces of the form $x+ \ker C$.}
\label{fig:effective_vector_field}
\end{figure}

In Figure~\ref{fig:effective_vector_field} we illustrate this by plotting the vector field 
\begin{align*}
x \mapsto \begin{pmatrix} 0 & 0 \\ 0 & 1 \end{pmatrix} Ax
\end{align*}
which is the ``effective component" of the drift vector field $\dot{x} = Ax$ in terms of change in output for this example.

This shows geometrically and very intuitively, why in our effort to reduce the value of the output by appropriately applying an input $u$, information of the state beyond its output value is completely irrelevant.

More abstractly, the invariance of $\ker C$ under the flow of $\dot{x}=Ax$ allows one to consider a well-defined flow on the quotient space $\mathbb R^n / \ker C$, cf. \cite{trentelman2012control},\cite{wonham2012linear}. To this end, let us consider the equivalence relation 
\begin{align*}
x' \sim x'' :\Leftrightarrow x' - x''\in \ker C,
\end{align*}
which is in fact the indistinguishability relation, since under the invariance of $\ker C$ under $\dot{x}=Ax$, the unobservable subspace is $\ker C$ itself (take $x_0 \ne 0$ in the kernel of $C$, then $ y(t)  := Cx(t) \equiv 0$). Moreover, $\mathbb R^n / \ker C$ is the state space $\mathbb R^n$ factored into indistinguishable equivalence classes. As we have seen, the effective velocity field in the orthogonal direction  to $\ker C$ (the direction that affects a change in the output)  for any points $x' \sim x''$ is equal. This is precisely the important fact that allows us to introduce a well-defined flow on the quotient $\mathbb R^n / \ker C$, which is in fact precisely \emph{represented} by the closed dynamics $\dot{y} = \tilde{A}y + CBu$. This gives a more complete view of the approach. Lastly, we note that the general idea presented in this section was also formulated in \cite{aoki1968control} in the context of model aggregation.

\section{Necessity of broadcast feedback in optimal centroid stabilization}
\label{sec:centroid}

In this section we apply the same principle as before to the optimal centroid stabilization problem. This will show that a broadcast feedback is necessarily the optimal structure. Note also that the optimal centroid stabilization problem can be viewed as the optimal stabilization of the (sample) mean of the group of linear systems, while the optimal synchronization problem presented in Section~\ref{sec:problem_formulation} can be viewed as the optimal stabilization of the (sample) variance.

We consider again a group of $N$ homogeneous linear systems
$$
 \dot{x}_i = A x_i + B u_i
$$
where $A \in \mathbb R^{n \times n}$, $B \in \mathbb R^{n \times p}$, and $(A,B)$ is controllable. Furthermore, we introduce the stacked system
$$
   \dot{x} = (I_N \otimes A) x + (I_N \otimes B) u.
$$
The goal in the optimal centroid stabilization problem is to stabilize in an optimal way the centroid of the group which is given by the arithmetic mean
$
  \bar{x} = \frac{1}{N} \sum_{i=1}^N x_i,
$
and can also be thought of as the center of mass. It is natural to introduce the fictitious output matrix
$
   C =  \begin{pmatrix} I & \dots & I \end{pmatrix}.
$
Following the static output feedback approach of the foregoing section, it is not difficult to see that
\begin{align*}
\dot{y} = A y + \begin{pmatrix} B & \dots & B \end{pmatrix} u.
\end{align*}
Therefore, the resulting optimal centroid stabilization problem can be cast as
\begin{align}
\begin{split}
    &\underset{u(\cdot)}{\text{minimize }}  \;  \int_{0}^{\infty}   (y^{\top} Q y + u^{\top} R u ) \, \text{d}t    \\
    &\text{subject to }  \;\; \dot{y} = A y +  \begin{pmatrix} B & \dots & B \end{pmatrix} u.
\end{split}
\label{eq:lqr_broadcast}
\end{align}
By the discussion of the general approach in the foregoing section, we can readily conclude that for $Q, R > 0$, the optimal feedback is of the form $u=Ky = K \bar{x}$. This feedback requires only aggregated information of the group, but the feedback is not quite a broadcast feedback yet, as at this point we only know that
\begin{align*}
 K = \begin{pmatrix} K_1 \\ \vdots \\ K_N \end{pmatrix},
\end{align*}
where it could be $K_i \ne K_j$. Note also that this already gives a (block) rank-one feedback matrix as advocated by \cite{madjidian2014distributed}. 
In the following we show that if the input weight $R$ is chosen homogeneously, then the resulting feedback is a broadcast feedback. More precisely, let $Q > 0$ and $R = (I_N \otimes W)$ with $W>0$. Then the structure of $R$ can be exploited in the solution of the LQR problem \eqref{eq:lqr_broadcast} given by
\begin{align*}
    u^{\star} = -R^{-1}  \begin{pmatrix} B & \dots & B \end{pmatrix}^{\top} Py,
\end{align*}
where $P$ denotes the unique positive definite solution of the algebraic Riccati equation associated to \eqref{eq:lqr_broadcast}.
Inserting $R^{-1} = (I_N \otimes W^{-1})$ one has
\begin{align*}
   u^{\star} = -  (\mathbf{1}_N \otimes W^{-1} B^{\top} P ) y.
\end{align*}
This finally yields the claimed broadcast structure. The resulting feedback can be rewritten by inserting $y = \begin{pmatrix} I & \dots & I \end{pmatrix} x$, which also showcases the (block) rank-one structure very clearly.

Furthermore, in the case of homogeneous input weights, the solution can be obtained very efficiently since the resulting algebraic Riccati equation for $P$ is given by
\begin{align*}
 PA + A^{\top} P - P \begin{pmatrix} B & \dots & B \end{pmatrix} \begin{pmatrix} W^{-1} B^{\top} \\ \vdots \\ W^{-1} B^{\top} \end{pmatrix} P + Q = 0.
\end{align*}
Therefore we only have to solve one ARE given by
\begin{align*}
PA + A^{\top} P - N P BW^{-1}B^{\top} P + Q = 0,
\end{align*}
which is of the dimension of \emph{one} individual system. In particular, the actual computation of the optimal control is scalable for the case of homogeneous weights.

One can think of the broadcast mechanism as having one central controller that exercises a global force that acts on all agents the same way. In particular, steering with a broadcast signal can be considered a coordination task for one central ``leader" (cf. \cite{brockett2010control}), which does not seem easy from an intuitive point of view. Nevertheless, our results show that in the case of centroid stabilization, a broadcast feedback naturally emerges as the optimal feedback structure, which is quite remarkable.

\section{Conclusions}
\label{sec:conclusions}
We considered two classes of LQR problems in distributed control, namely optimal synchronization and optimal center of mass stabilization. By employing an approach that translates the problems into a special LQR problem in which the cost functional integrates over squared outputs and inputs, and exploiting the structure of the problems, we were able to give a unifying solution based on showing that the optimal control is always a static output feedback. For the case of optimal synchronization with homogeneous agents we showed that the optimal solution is given by a scalable diffusive coupling. For the case of optimal center of mass stabilization, we obtained as the optimal solution a broadcast feedback which has nice scalable properties as well. The general LQR problem that considers outputs in the cost functional as well as our approach employing this viewpoint is also of interest in its own right. We formulated the general underlying principle and also discussed and illuminated it from a geometric point of view. An interesting point for future research is the study of this underlying principle in more general cases that involve e.g. nonlinear systems.

\section{Acknowledgments}
We are indebted to Jan Maximilian Montenbruck for the very interesting and fruitful discussions on this work.

\appendix
\section{Linear quadratic regulator problems}
In this subsection, we gather relevant results from linear quadratic regulator theory, see also the appendix in \cite{Montenbruck}. We consider a linear system
$
\dot{x} = Ax +  Bu
$
and with $Q \ge 0$ and $R>0$, a cost functional
$$
  J = \int_{0}^{\infty} (x^{\top}Qx + u^{\top} R u ) \, \text{d}t.
$$
The linear quadratic regulator problem seeks the input $u$ which minimizes the cost functional. A fundamental result in linear quadratic regulator theory is the connection between the optimal solution $u$ and the algebraic Riccati equation
\begin{align}
   XBR^{-1}B^{\top}X - XA - A^{\top}X - Q = 0.
\label{eq:ARE_general}
\end{align}
A very general result in linear quadratic regulator theory relates the optimal solution $u$ to the so-called \emph{smallest positive semidefinite solution} of the algebraic Riccati equation, which is a solution $X_{-}$ of \eqref{eq:ARE_general} which satisfies $X_{-} \le \tilde{X}$ for any other solution $\tilde{X}$ of \eqref{eq:ARE_general}. Given this smallest positive semidefinite solution $X_{-}$ of the algebraic Riccati equation \eqref{eq:ARE_general}, the control
\begin{align}
    u^{\star} = -R^{-1} B^{\top} X_{-} x
 \label{eq:optimal_control_smallest_positive}
\end{align}
minimizes the cost functional $J$, see e.g.\ \cite{trentelman1989regular}, Theorem 4.2.

Furthermore the notion of strong and maximal solutions is useful. A solution $X_s$ to the algebraic Riccati equation \eqref{eq:ARE_general} is called \emph{strong}, if $A - BR^{-1}B^{\top} X_s$ has all of its eigenvalues in the closed left half-plane. A symmetric solution $X_{+}$ to the algebraic Riccati equation \eqref{eq:ARE_general} is said to be \emph{maximal}, if $X_{+} > \tilde{X}$ for every other symmetric solution $\tilde{X}$.

Let $G G^{\top} = R^{-1}$ with $G$ full rank and $F^{\top} F = Q$ with $F$ full rank. The following results establish a precise connection between the three different types of solutions of the algebraic Riccati equation.

\begin{Theorem}[\cite{molinari1977time}] 
Let  $(A,BG)$ be stabilizable. If the algebraic Riccati equation \eqref{eq:ARE_general} has a strong solution $X_s$, then $X_s$ is a maximal solution. If $(A,BG)$ is stabilizable, then \eqref{eq:ARE_general} has at most one strong solution.
\end{Theorem}

\begin{Theorem}[\cite{wimmer1985monotonicity}]
Let $(A,B)$ be stabilizable and let $X_{-}$ and $X_{+}$ denote the smallest positive semidefinite and maximal solution to the algebraic Riccati equation \eqref{eq:ARE_general}, respectively. Then $X_{-} = X_{+}$ if and only if the space spanned by $(A,F)$-undetectable eigenvectors associated with eigenvalues of $A$ in the open right half-plane is trivial.
\label{thm:all_three_solutions}
\end{Theorem}

Thus, if, for instance, $A$ has only stable eigenvalues, then the condition of the above theorem is naturally fulfilled so that $X_{-}$ in \eqref{eq:optimal_control_smallest_positive} can be replaced with the maximal solution $X_{+}$ which in turn is equal to the strong solution $X_s$ if it exists. Furthermore if $Q>0$, there are no $(A,F)$-undetectable eigenvectors and the aforementioned holds true as well.

\section{Proof of Proposition~\ref{prop:scalability}}
\label{appendix:proof}
In the case of homogeneous weights, the algebraic Riccati equation for the considered optimal control problem reads as
\begin{align*}
 P (\Gamma^{\top} \otimes B)& (I_{N-1} \otimes W)^{-1} (\Gamma^{\top} \otimes B)^{\top} P  \\
  &\hspace{-1.6cm} - P (I_{N-1} \otimes A) - (I_{N-1} \otimes A)^{\top} P - (I_{N-1} \otimes V) = 0.
\end{align*}
Inverting $(I_{N-1} \otimes W)^{-1} = I_{N-1} \otimes W^{-1} $ and using the mixed-product property of the Kronecker product, as well as $\Gamma^{\top}\Gamma = I_{N-1}$, we rewrite the ARE as
\begin{align}
\begin{split}
P( I_{N-1}  &\otimes BW^{-1}B^{\top})P - P(I_{N-1} \otimes A)  \\
&\hspace{0.5cm}  - (I_{N-1} \otimes A^{\top}) - P (I_{N-1} \otimes V) = 0.
\end{split}
\label{eq:ARE_decoupled}
\end{align}
The block diagonal structure of the matrices on the left-hand side results in the fact that \eqref{eq:ARE_decoupled} is decoupled. Furthermore, the individual AREs are all the same. This suggests the ansatz $P = I_{N-1} \otimes Y$ which, inserted in \eqref{eq:ARE_decoupled}, leads to
\begin{align}
  Y BW^{-1}B^{\top} Y - YA - A^{\top} Y - V = 0.
\label{eq:reduced_ARE}
\end{align}
This is the algebraic Riccati equation associated to a different (fictitious) LQR problem with weights $V, W$ and dynamics $(A,B)$. Since this LQR problem fulfills all the assumptions of the standard LQR setup, we conclude that \eqref{eq:reduced_ARE} has a unique solution $Y^{\star}$, and hence $P^{\star} = I_{N-1} \otimes Y^{\star}$ is the unique positive definite solution of \eqref{eq:ARE_decoupled}. Moreover, applying again the mixed product rule to $K = - R^{-1}(\Gamma^{\top} \otimes B)P^{\star}$, we get
$
 K   = -   (\Gamma \otimes W^{-1} B^{\top} Y^{\star}),
$
and thus
\begin{align*}
     u =  -  K (\Gamma^{\top} \otimes I) x
        = - (\Gamma \Gamma^{\top} \otimes W^{-1} B^{\top} Y^{\star}) x.
\end{align*}
Due to the definition of $\Gamma$, the matrix $\Gamma \Gamma^{\top}$ is the orthogonal projection $P_{\mathscr A'}$ on the asynchronous subspace for $N$ scalar systems.
Since $P_{\mathscr A'} = I_N - P_{\mathscr S'}$ and 
$
P_{\mathscr S'} = \big( \frac{1}{\sqrt{N}} \mathbf{1}_N \big) \big(\frac{1}{\sqrt{N}} \mathbf{1}_N\big)^{\top},
$
we see that $P_{\mathscr A'} = \frac{1}{N} \mathcal L$, where $\mathcal L$ denotes the graph Laplacian of the complete graph.

To conclude, in the case of homogeneous weights, the dimension of the ARE that needs to be considered is reduced from $Nn$ to $n$, see also \cite{Montenbruck}. Furthermore we can conclude that the resulting feedback structure can be written as
\begin{align*}
  u_i = \frac{1}{N} \sum_{j=1}^N   W^{-1} B^{\top} Y^{\star} (x_j - x_i),
\end{align*}
i.e. the diffusive law has the further property that the coupling gains between two coupled systems are identical.  \hfill \phantom{.}  \qed


\bibliographystyle{ieeetr}
\bibliography{promotion_library}

\end{document}